\documentclass[12pt,fullpage]{article}

\usepackage{amsfonts}
\usepackage{xr}
\usepackage{graphicx}
\usepackage{amsmath}
\usepackage{epsfig,subfigure}
\usepackage{color}
\def\R{\mathbb{R}}
\def\N{\mathbb{N}}

\def\epsilon{\varepsilon}
\def\hat{\widehat}
\def\tilde{\widetilde}

\newcommand{\SE}{\setcounter{equation}{0} \section}
\newcommand{\be}{\begin{equation}}
\newcommand{\ee}{\end{equation}}
\newcommand{\baa}{\begin{array}}
\newcommand{\eaa}{\end{array}}
\newcommand{\ba}{\begin{eqnarray}}
\newcommand{\ea}{\end{eqnarray}}

\newtheorem{theo}{\bf Theorem}[section]
\newtheorem{lem}[theo]{\bf Lemma}

\newtheorem{rem}[theo]{\bf Remark}

\begin{document}
\date{}
\title{\bf{On the nonlocal Fisher-KPP equation: steady states,
    spreading speed and 
global bounds}}
\author{Fran\c cois Hamel\thanks{Aix Marseille Universit\'e, CNRS, Centrale Marseille, LATP, UMR 7353, 13453 Marseille, France, and Institut Universitaire de France}   
\and Lenya Ryzhik\thanks{Department of Mathematics, Stanford University, Stanford, CA 94305, USA}}
 
\maketitle

\begin{abstract}
We consider the Fisher-KPP equation with a non-local interaction
term. We establish a condition on the interaction that allows for
existence of non-constant periodic solutions, and prove uniform upper
bounds for the solutions of the Cauchy problem, as well as upper and
lower bounds on the spreading rate of the solutions with compactly
supported initial data. 
\end{abstract}


\SE{Introduction and main results}\label{intro}

Reaction-diffusions equations of the form
\be\label{eq}
u_t=u_{xx}+\mu\,u\,(1-\phi*u),\ \ t>0,\ \ x\in\R
\ee
model the non-local interaction and competition between species. Here,
$u(t,x)$ is a population density, and the independent variable $x$ may
be either a spatial location, a trait, or a mutation of the
species. This equation is a generalization of the classical (local)
Fisher-KPP (for Kolmogorov-Petrovsky-Piskunov) \cite{f,kpp} equation
\begin{equation}\label{fisher-local}
u_t=u_{xx}+\mu u(1-u).
\end{equation}
The parameter $\mu>0$ in both cases denotes the strength of the
competition. The nonnegative convolution kernel $\phi\in L^1(\R)$ in (\ref{eq})
models the species interaction. We assume that it satisfies the
following properties:
\begin{equation}\label{hypphi0}
\hbox{$\phi(x)\ge 0$ for all $x\in\R$}, ~~~
\mathop{{\rm{ess}}\,{\rm{inf}}}_{(-\sigma,\sigma)}\phi>0
\hbox{ for some }\sigma>0,\ \hbox{ and }\ \int_{\R}\phi=1.
\end{equation}
This equation arises, in particular, in ecology, and adaptive
dynamics~\cite{FG,GVA,Gourley}, but also in essentially any area where
the local Fisher-KPP equation (\ref{fisher-local}) appears, as
interactions are rarely fully local, and often occur on scales that
are comparable to the diffusion scale, making the non-local equation a
more realistic model.

\subsubsection*{Qualitative behavior of traveling waves and numerical solutions}

While the behavior of the solutions of the local equation
(\ref{fisher-local}) has been very extensively studied, much less is
known about solutions of (\ref{eq}). Numerical
simulations~\cite{bnpr,GVA,Gourley} indicate the following behavior of
solutions: first, for small $\mu$ solutions of the non-local equation
behave ``exactly like solutions of the local equation''. More
precisely, numerics shows that when $\mu$ is small, the only
non-negative bounded steady states for (\ref{eq}) are $u\equiv 0$ and
$u\equiv 1$, and this equation admits monotonic traveling waves of the
form
\[
u(t,x)=U(x-ct),\ \ U(-\infty)=1,\ \ U(+\infty)=0,\ \ U\hbox{ is decreasing},
\]
that are asymptotically stable for the solutions of the Cauchy problem
with front-like initial data. On the other hand, when $\mu$ is
sufficiently large, the qualitative behavior of numerical solutions
for (\ref{eq}) and (\ref{fisher-local}) is very different: first,
non-constant bounded steady solutions of (\ref{eq}) may exist. As far
as traveling waves are concerned, when $\mu$ is large, numerical
computations show that non-monotonic traveling waves connecting $1$ to
$0$ may exist, as well as non-monotonic pulsating waves
\[
u(t,x)=U(x-ct,x),\ \ U(-\infty,x)>0,\ \ U(+\infty,x)=0,\ \ U\hbox{ is periodic in }x,
\]
connecting non-constant steady states to $u\equiv 0$. Moreover, the
long time limit of the solutions of the Cauchy problem may be either
non-monotonic traveling waves, pulsating waves or wave-trains
\[
u(t,x)=U(x-ct),\ \ U\hbox{ is periodic},
\]
and not the monotonic traveling waves (see~\cite{bnpr,GVA,Gourley},
and also~\cite{abvv} for a related model). To the best of our
knowledge, so far, these numerical
observations have not been proved rigorously. 

When $\mu$ is sufficiently large or when $\phi$ is sufficiently far
from a Dirac mass at $0$, the steady solution $u\equiv 1$ becomes
unstable with respect to some spatially periodic perturbations,
provided that the Fourier transform $\hat\phi(\xi)$ is not
non-negative everywhere. This, heuristically, leads to existence of
the non-constant periodic steady states or wave trains that were
observed numerically~\cite{GVA,Gourley} (see also~\cite{b,gb} for a
bifurcation analysis and numerics on a related equation). Quite surprisingly, it
was shown numerically in~\cite{NPT} that even when the state $u\equiv
1$ is unstable, there exist traveling wave solutions (which are
themselves not stable) that connect $u\equiv 1$ to $u\equiv 0$. It has
been proved in~\cite{bnpr} that when $\hat\phi(\xi)\ge 0$ for all
$\xi\in\R$, so that $u\equiv 1$ is stable, and $\mu$ is sufficiently
large, monotonic traveling waves connecting $u\equiv 0$ and $u\equiv
1$ do not exist, but there exists a non-monotonic traveling wave
connecting these two states. We also point out that explicit examples
of various wave-trains have been recently constructed in~\cite{NPRR}.

As far as monotone traveling waves for (\ref{eq}) are concerned, Fang
and Zhao \cite{FangZhao} have proved that for all $c\geq 2$, there
exists $\mu_c>0$, so that (\ref{eq}) admits a monotonic traveling wave
if and only if $\mu \in (0,\mu_c)$ -- see also~\cite{Gourley,wlr}.  It has been also shown in~\cite{AlfaroCoville} that
traveling waves with large speeds necessarily converge to $1$ as $x-ct\to-\infty$.  

Essentially nothing is known about the long time behavior of solutions
of the Cauchy problem for (\ref{eq}), in particular, whether the
aforementioned traveling waves are stable.  The main mathematical
difficulty in studying this question is that unlike for
(\ref{fisher-local}), solutions of (\ref{eq}) do not obey the maximum
principle or the comparison principle, making the use of many
technical tools impossible.   

\subsubsection*{Non-trivial steady states}

We study here existence of nonnegative steady bounded solutions $u=u(x)$ of
(\ref{eq}):
\begin{equation}\label{stat}
u''(x)+\mu\,u(x)\,(1-(\phi*u)(x))=0,\ \ x\in\R.
\end{equation}
Solutions are understood in the classical sense $C^2(\R)$. Notice
that, for any such solution $u$,~(\ref{stat}) can be written as
\[
\hbox{$u''+cu=0$ in $\R$}, 
\]
where $c=\mu\,(1-\phi*u)$ is continuous and
bounded in $\R$. Hence, for any nonnegative bounded solution $u$
of~(\ref{stat}), either $u=0$ in $\R$ or $u>0$ in $\R$, from the
strong maximum principle.
The constant functions $u=0$ and $u=1$ solve~(\ref{stat}) and they are
the only constant functions solving~(\ref{stat}), for every
$\mu>0$. The following theorem gives a sufficient condition for other
non-trivial periodic steady states to exist when $\mu$ is sufficiently
large.  
We denote by $\hat{\phi}$ the Fourier transform of the function
$\phi$:
\[
\hat{\phi}(\xi)=\int_{\R}\phi(x)e^{-2i\pi\xi x}dx\ \hbox{ for all
}\xi\in\R.
\]
\begin{theo}\label{th1}
Assume that $\phi$ is even and that there 
are $L\in(0,+\infty)$ and $k_0\in\N$ such that
\begin{equation}\label{hypphi}
\hat{\phi}\Big(\frac{k_0}{L}\Big)<0\ \hbox{ and }\ 
\hat{\phi}\Big(\frac{k}{L}\Big)\ge0\hbox{ for all }k\in\N\backslash\{k_0\}.
\end{equation}
Then there is $\mu^*>0$ such that, for every $\mu>\mu^*$, the
problem~$(\ref{stat})$ has a non-constant positive $L$-periodic
solution~$u$.
\end{theo}
A sufficient condition for~(\ref{hypphi}) is that there exists
$\xi^*>0$ such that $\hat{\phi}\ge0$ on $[\xi^*,+\infty)$
and~$\hat{\phi}<0$ in a left neighborhood of $\xi^*$. In that case, we
can take $L>0$ such that $1/L<\xi^*<2/L$ and $\hat{\phi}(1/L)<0$, that
is~(\ref{hypphi}) holds with $k_0=1$. For instance, consider $\beta>1$
large enough so that 
\[
c_{\beta}:=1-1/\sqrt{\beta}+1/\beta>0,
\]
and define
\[
\phi_{\beta}(x)=\frac{1}{c_{\beta}\sqrt{\pi}}\big(e^{-x^2}
-e^{-\beta x^2}+e^{-\beta^2x^2}\big).
\]
The function $\phi_{\beta}$ is positive in $\R$ and its $L^1(\R)$ norm is equal to $1$. Furthermore, the function~$\hat{\phi}_{\beta}$ is such that
$$\hat{\phi}_{\beta}(\sqrt{\beta})=\frac{1}{c_{\beta}}\Big(e^{-\beta\pi^2}-\frac{e^{-\pi^2}}{\sqrt{\beta}}+\frac{e^{-\pi^2/\beta}}{\beta}\Big)\sim-\frac{e^{-\pi^2}}{\sqrt{\beta}}<0\ \hbox{ as }\beta\to+\infty,$$
while $\hat{\phi}_{\beta}(\xi)\sim(c_{\beta}\beta)^{-1}e^{-\pi^2\xi^2/\beta^2}>0$ as $\xi\to+\infty$ for every fixed $\beta>1$ such that $c_{\beta}>0$. As a consequence, the function $\phi_{\beta}$ satisfies the assumptions of Theorem~\ref{th1} for $\beta>1$ large enough.

\subsubsection*{Bounds for the solutions of the Cauchy problem}

We also consider solutions of the  
Cauchy problem
\begin{equation}\label{cauchy}\left\{\baa{l}
u_t=u_{xx}+\mu\,u\,(1-\phi*u),\ \ t>0,\ \ x\in\R,\vspace{3pt}\\
u(0,\cdot)=u_0,\eaa\right.
\end{equation}
with a non-negative initial condition $u_0\in L^{\infty}(\R)$, $u_0\ge
0$. Solutions of the corresponding Cauchy problem for the local
Fisher-KPP equation~(\ref{fisher-local}) satisfy a trivial upper bound 
\[
u(t,x)\le \max(1,\|u_0\|_{L^\infty}).
\]
However, as the maximum principle does not hold for (\ref{cauchy}),
such bound needs not hold for the non-local problem. We prove here the
following result. 
\begin{theo}\label{th2}
  For every $\mu>0$ and for every nonnegative initial
  condition~$u_0\in L^{\infty}(\R)$, the solution $u$
  of~$(\ref{cauchy})$ exists and is globally bounded in time, that is
  there exists $M>0$ such that
\[
0\le u(t,x)\le M\ \hbox{ for all }t>0\hbox{ and }x\in\R.
\]
\end{theo}
In Theorem~$\ref{th2}$, the condition that $\phi\ge0$ belongs to
$L^1(\R)$ with unit mass could likely be replaced by the assumption that
$\phi$ is a nonnegative Radon probability measure. On the other hand,
the condition 
\[
{\rm{ess}}\,{\rm{inf}}_{(-\sigma,\sigma)}\phi>0
\] 
for some $\sigma>0$ plays a crucial role in
Theorem~$\ref{th2}$. Without it, the conclusion fails in general. Let
us illustrate this with the following example. Let $\mu$ and $L$ be two
positive real numbers such that $\mu>\pi^2/L^2$ and consider
(\ref{cauchy}) with
\[
\phi=\frac{1}{2}\big(\delta_{-L}+\delta_L\big),
\]
where $\delta_{\pm L}$ denote the Dirac masses at the points $\pm
L$. Assume that the initial data $u_0$ is $2L$-periodic. By
uniqueness, the function $u(t,\cdot)$ is $2L$-periodic for every $t>0$
and it obeys
\begin{eqnarray*}
&&u_t(t,x) = u_{xx}(t,x)+\mu\,u(t,x)\,
\Big(1-\displaystyle\frac{u(t,x-L)+u(t,x+L)}{2}\Big)\vspace{3pt}\\
&&~~~~~~~~~ = u_{xx}(t,x)+\mu\,u(t,x)\,(1-u(t,x+L)),
\end{eqnarray*}
for every $(t,x)\in(0,+\infty)\times\R$. Set $v(t,x)=u(t,x+L)$ and
$w(t,x)=u(t,x)-v(t,x)$. The functions $u$, $v$ and $w$ satisfy
\[
\baa{rcl}
&&u_t  =  u_{xx}+\mu\,u\,(1-v),\vspace{3pt}\\
&&v_t  =  v_{xx}+\mu\,v\,(1-u),\vspace{3pt}\\
&&w_t  =  w_{xx}+\mu\,w,\eaa
\]
in $(0,+\infty)\times\R$. Therefore, if for instance one chooses
\[
u_0(x)=1+\rho\,\cos\Big(\frac{\pi x}{L}\Big)
\]
for some $0<\rho<1$, then $w(0,x)=2\,\rho\,\cos(\pi x/L)$ and
\[
w(t,x)=2\,\rho\,\cos\Big(\frac{\pi x}{L}\Big)\,e^{(\mu-\pi^2/L^2)t},
\]
for every $t>0$ and $x\in\R$. Since $\mu>\pi^2/L^2$, one gets 
\[
\hbox{$\|w(t,\cdot)\|_{L^{\infty}(\R)}\to+\infty$ as $t\to+\infty$.}
\] 
However, we have 
\[
\|u(t,\cdot)\|_{L^{\infty}(\R)}=\|v(t,\cdot)\|_{L^{\infty}(\R)},
\]
and 
\[
\|w(t,\cdot)\|_{L^{\infty}(\R)}\le\|u(t,\cdot)\|_{L^{\infty}(\R)}+\|v(t,\cdot)\|_{L^{\infty}(\R)}=2\|u(t,\cdot)\|_{L^{\infty}(\R)},
\] 
for every $t>0$. Finally, one concludes that $\|u(t,\cdot)\|_{L^{\infty}(\R)}\to+\infty$ as $t\to+\infty$.

\subsubsection*{Spreading rate for compactly supported initial conditions}

Finally, we establish the following result for the spreading rate for
solutions of the Cauchy problem with compactly supported initial data.
\begin{theo}\label{th3}
Let $u$ be the solution of the Cauchy problem~$(\ref{cauchy})$ with a nonnegative initial condition $u_0\in L^{\infty}(\R)$ such that $u_0\not\equiv 0$. Then
\begin{equation}\label{spreadinf}
  \liminf_{t\to+\infty}\Big(\min_{|x|\le ct}u(t,x)\Big)>0\ \hbox{ for all }0\le c<2\sqrt{\mu}.
\end{equation}
Furthermore, if $u_0$ is compactly supported, then
\begin{equation}\label{spreadsup}
\lim_{t\to+\infty}\Big(\max_{|x|\ge ct}u(t,x)\Big)=0\ \hbox{ for all }c>2\sqrt{\mu}.
\end{equation}
\end{theo}
Similar bounds for the spreading rate for the solutions of the local
Fisher-KPP equation (\ref{fisher-local}) are well known~\cite{aw}. And, indeed,  
the upper bound (\ref{spreadsup}) on the spreading rate  follows
immediately from the comparison of the solutions of (\ref{cauchy}) and
the solutions of the linear problem 
\begin{equation}\label{linear-local}
v_t=v_{xx}+\mu v.
\end{equation}
On the other hand, the lack of the maximum principle for
(\ref{cauchy}) makes the proof of the lower bound~(\ref{spreadinf})
very different from that for the local equation
(\ref{fisher-local}). The idea is, roughly, as follows. If $u(t,x)$ is
small behind the location $x(t)=2\sqrt{\mu}t$ then it should be small on
a sufficiently large region to make also the convolution $\phi*u$
small.  In that case, however, $u(t,x)$ behaves, approximately, as a
solution of the linearized equation (\ref{linear-local}). These
solutions, however, propagate with the speed $c_*=2\sqrt{\mu}$,
leading to a contradiction.

The paper is organized as follows. Section~\ref{sec2} contains the
proof of Theorem~\ref{th1}, while Theorems~\ref{th2} and~\ref{th3} are
proved in Sections~\ref{sec3} and~\ref{sec4}, respectively.\hfill\break

\noindent{\bf Acknowledgment.} The research leading to these results has
received funding from the French ANR within the project PREFERED and
from the European Research Council under the European Union's Seventh
Framework Programme (FP/2007-2013) / ERC Grant Agreement n.321186~-~ReaDi~-~Reaction-Diffusion Equations, Propagation and Modelling. Part
of this work was also carried out during visits by F.~Hamel to the
Departments of Mathematics of the University of California, Berkeley
and of Stanford University, the hospitality of which is thankfully
acknowledged.  This work was also supported by the NSF grants
DMS-0908507 and DMS-115893. The authors are grateful to G.~Nadin and B.~Perthame for valuable discussions during the preparation of this work.


\SE{Existence of non-constant periodic steady states}\label{sec2}

This section contains the proof of Theorem~\ref{th1}. The general strategy is to apply a topological degree argument. We first prove that, if $\mu>0$ is small enough, the constant solution $u\equiv 1$ is linearly stable and is the only bounded positive solution of~(\ref{stat}). Next, we show that, under the assumptions of the theorem, if $\mu>0$ is large enough, the constant~$u\equiv 1$ is an isolated (in the $L^{\infty}$-sense) unstable steady solution with one unstable direction in the class of $L$-periodic even functions. Together with some a priori uniform estimates, this will allow us conclude that $u\equiv 1$ is not the only positive solution of (\ref{stat}) for large values of the parameter $\mu$. Throughout this section, we assume that the convolution kernel $\phi$ is a nonnegative $L^1(\R)$ function satisfying~(\ref{hypphi0}). The additional assumption (\ref{hypphi}) and the evenness of $\phi$ will be used only in some steps of the proof.

\subsubsection*{Uniform bounds for solutions}

We first establish uniform pointwise upper and lower bounds for all positive bounded solutions of~(\ref{stat}), when the parameter $\mu$ ranges between two fixed positive constants.
\begin{lem}\label{lem1}
For every $0<a\le b<+\infty$, there is a positive constant $m$ which depends on~$\phi$,~$a$ and~$b$ and there is a positive constant $M$ which depends on~$\phi$ and~$b$, such that, for every $\mu\in[a,b]$ and every positive bounded solution $u$ of~$(\ref{stat})$, there holds
\be\label{mM}
0<m\le u(x)\le M\ \hbox{ for all }x\in\R.
\ee
\end{lem}

\noindent{\bf{Proof.}} Similar lower and upper bounds were proved for a fixed value of $\mu>0$ in Lemma~2.1 and Proposition~2.3 of~\cite{bnpr}, under the additional assumption that the convolution kernel~$\phi$ is of class~$C^1$ and has a bounded second moment. Lemma~\ref{lem1} of the present paper can be viewed as an extension of the aforementioned results. Furthermore, the proof of the upper bound is slightly different from that of Proposition~2.3 of~\cite{bnpr}.

We begin with the proof of the upper bound, that is the existence of a constant $M$ in~(\ref{mM}), independently of $\mu\in[a,b]$ and of $u>0$ solving~(\ref{stat}). Let $\sigma>0$ be as in~(\ref{hypphi0}), and $\delta$ be any fixed real number such that 
\begin{equation}\label{delta} 0<\delta<\min\Big(\frac{\pi}{2\sqrt{b}},\sigma\Big).
\end{equation}
Finally, let $\eta>0$ be such that 
\begin{equation}\label{eta} 
\phi(x)\ge\eta>0\ \hbox{ for almost every }x\in(-\delta,\delta).  
\end{equation} 
From the choice of $\delta$, we know that the lowest eigenvalue of the operator 
\begin{equation}\label{oper-max}
\varphi\mapsto-\varphi''-b\varphi
\end{equation}
in $(-\delta,\delta)$ with the Dirichlet boundary condition at $\pm\delta$, which is 
equal to $\pi^2/(4\delta^2)-b$, is positive. Therefore,   the weak maximum principle for 
this operator holds in the interval $(-\delta,\delta)$ (see \cite{bnv}) and in any subinterval, 
and there is a unique $C^2([-\delta,\delta])$ solution $\overline{u}$ of the boundary value 
problem
\begin{equation}\label{baru-eq}
\overline{u}''+b\overline{u}=0\hbox{ in }[-\delta,\delta]\ \hbox{ and }\ \overline{u}(\pm\delta)=\frac{1}{\delta\eta}.
\end{equation}
Furthermore, since the constant $1/(\delta\eta)$ is a subsolution for the above equation, 
the weak maximum principle yields
\begin{equation}\label{baru-bd}
\overline{u}\ge\frac{1}{\delta\eta}\hbox{ in }[-\delta,\delta].
\end{equation}
Notice that $\delta$ depends only on $\phi$ and $b$, that $\eta$ depends only on 
$\phi$ and $\delta$ and that $\overline{u}$ depends only on~$b$,~$\delta$ and~$\eta$, 
that is $\overline{u}$ depends only on $\phi$ and $b$.

Let now $\mu$ be any real number in $[a,b]$, $u$ be any positive 
bounded solution of~(\ref{stat}) and set
\[
\rho_u=\sup_{\R}u>0.
\]
Since (\ref{stat}) implies that $u''$ is bounded, we know that $u'$ is also bounded.
By differentiating~(\ref{stat}), we see that the $k$th-order derivative $u^{(k)}$ of 
$u$ exists and is 
bounded, for every $k\in\N$. Furthermore, there is a sequence $(x_n)_{n\in\N}$ of real numbers 
such that
\[
u(x_n)\to\rho_u\hbox{ as }n\to+\infty.
\]
Consider the translates
\[
u_n(x)=u(x+x_n)\ \hbox{ for }n\in\N\hbox{ and }x\in\R.
\]
Up to extraction of a subsequence, the functions $u_n$ converge in (at least) 
$C^2_{loc}(\R)$ to a $C^{\infty}(\R)$ solution $u_{\infty}$ of~(\ref{stat}) such that 
$0\le u_{\infty}\le\rho_u$ in $\R$ and $u_{\infty}(0)=\rho_u>0$. Therefore, 
we have $u_{\infty}''(0)\le0$, whence
\[
\int_{\R}\phi(y)\,u_{\infty}(-y)\,dy=(\phi*u_{\infty})(0)\le 1.
\]
Since both $\phi$ and $u_{\infty}$ are nonnegative, one gets that
\[
\int_{-\delta}^0\phi(y)\,u_{\infty}(-y)\,dy\le 1
\hbox{ and }\int_0^{\delta}\phi(y)\,u_{\infty}(-y)\,dy\le 1.
\]
Here, $\delta$ is defined in~(\ref{delta}). On the other hand, 
$\phi\ge\eta>0$ almost everywhere in $(-\delta,\delta)$, due to~(\ref{eta}). 
One infers the existence of some real numbers $y_{\pm}$ such that
\begin{equation}\label{u-interval}
-\delta\le y_-< 0< y_+\le\delta\ \hbox{ and }\ u_{\infty}(y_{\pm})\le\frac{1}{\delta\eta}.
\end{equation}
On the other hand, the nonnegative function $u_{\infty}$ satisfies
\[
0=u_{\infty}''+\mu\,u_{\infty}\,(1-\phi*u_{\infty})\le u_{\infty}''+\mu\,u_{\infty}\le u_{\infty}''+b\,u_{\infty}\hbox{ in }\R.
\]
That is, $u_\infty$ is a sub-solution to the equation (\ref{baru-eq}) satisfied by 
the function $\overline{u}$. Therefore, the lower bound (\ref{baru-bd}) on $\overline{u}$, 
and the upper bound (\ref{u-interval}) for $u_\infty$ at $y=y_{\pm}$, together 
with the maximum principle for the operator (\ref{oper-max}), imply that 
\[
u_{\infty}\le\overline{u}\hbox{ in }[y_-,y_+],
\]
whence $\rho_u=u_{\infty}(0)\le\overline{u}(0)$.  Since $\overline{u}(0)$ depends only on $\phi$ and $b$ and does not depend on $\mu$ nor on $u$, the upper bound in Lemma~\ref{lem1} is thereby proved.

We now turn to the proof of the lower bound in Lemma~\ref{lem1}, that is, the existence of 
a positive constant $m$, independent of $\mu\in[a,b]$ so that any bounded solution $u>0$ of 
(\ref{stat}) is bounded from below by $m$. Assume that there is no such constant $m>0$. 
Then, there exist a sequence $(\mu_n)_{n\in\N}$ of real numbers in the interval $[a,b]$, 
a sequence $(u_n)_{n\in\N}$ of positive bounded solutions of~(\ref{stat}) with 
$\mu=\mu_n$, and a sequence $(x_n)_{n\in\N}$ of real numbers such that
\[
u_n(x_n)\to0\hbox{ as }n\to+\infty.
\]
Each function $u_n$ is positive, and, from the already proved part of the present lemma, 
there is a constant $M>0$ such that $u_n\le M$ in $\R$ for all $n\in\N$. As a consequence, 
the sequence $(u_n'')_{n\in\N}$ is bounded in $L^{\infty}(\R)$, as are $(u_n')_{n\in\N}$ 
and $(u_n^{(k)})_{n\in\N}$ for every $k\in\N$ by immediate induction. Again, we consider the shifts
\[
v_n(x)=u_n(x+x_n).
\]
Up to extraction of a subsequence, one can assume that $\mu_n\to\mu\in[a,b]$ as $n\to+\infty$ and that the sequence $(v_n)_{n\in\N}$ converges as $n\to+\infty$ in (at least) $C^2_{loc}(\R)$ to a $C^{\infty}(\R)$ solution $v_{\infty}$ of~(\ref{stat}) such that $0\le v_{\infty}\le M$ in $\R$ and $v_{\infty}(0)=0$. The function $v_{\infty}$ satisfies 
\[
v_{\infty}''+cv_{\infty}=0\hbox{ in $\R$,}
\]
where
\[
c=\mu\,(1-\phi*v_{\infty})
\]
is a bounded continuous function. The strong maximum principle implies that $v_{\infty}\equiv 0$ 
in $\R$, that is~$v_n\to0$ as $n\to+\infty$ (at least) locally uniformly in $\R$. Since $\phi$ 
is a nonnegative $L^1(\R)$ function, and all functions $v_n$ satisfy $0<v_n\le M$ in $\R$, 
it follows that
\[
\phi*v_n\to0\hbox{ as }n\to+\infty\hbox{ locally uniformly in }\R.
\]
In particular, for every $R>0$, there is $N\in\N$ such that $\phi*v_n\le1/2$ in $[-R,R]$ 
for every $n\ge N$, whence
\[
0=v_n''+\mu_n\,v_n\,(1-\phi*v_n)\ge v_n''+\frac{\mu_nv_n}{2}\ge v_n''+\frac{av_n}{2}\ 
\hbox{ in }[-R,R]
\]
for every $n\ge N$. In other words, all functions $v_n$ (for $n\ge N$) are supersolutions of the elliptic operator 
\[
\varphi\mapsto-\varphi''-(a/2)\varphi ~~\hbox{ in $[-R,R]$}.
\]
Since the functions $v_n$ are all positive in $[-R,R]$, it follows from~\cite{bnv} that 
the lowest eigenvalue of this operator with Dirichlet boundary condition at the points~$\pm R$ 
is positive, that is
\[
\frac{\pi^2}{4R^2}-\frac{a}{2}>0.
\]
Since this holds for every $R>0$, one gets a contradiction as $R\to+\infty$. Hence, 
the sequences $(\mu_n)_{n\in\N}$, $(u_n)_{n\in\N}$ and $(x_n)_{n\in\N}$ as above can not exist. 
That yields the existence of $m$ in~(\ref{mM}) and completes the proof of Lemma~\ref{lem1}.\hfill$\Box$

\subsubsection*{Non-existence of non-trivial solutions for small $\mu$}

We proceed with the proof of Theorem~\ref{th1}.  Let now $L>0$ and $k_0\in\N$ be as in the assumption~(\ref{hypphi}) of Theorem~\ref{th1}.  We also assume from now on that the convolution kernel $\phi$ is even. Let $\alpha$ be any fixed real number in the interval~$(0,1)$, and let $X$ be the Banach space of even $L$-periodic functions of class~$C^{0,\alpha}(\R)$, equipped with the norm \[ 
\|u\|_X=\|u\|_{L^{\infty}(\R)}+\sup_{x\neq y\in\R}\frac{|u(x)-u(y)|}{|x-y|^{\alpha}}.  
\]

For every $\mu>0$, let $T_{\mu}:X\to X$ be an operator defined as follows: for every 
$u\in X$, $v:=T_{\mu}(u)\in X$ is the unique classical solution of
\[
-v''+v=u+\mu\,u\,(1-\phi*u)\hbox{ in }\R.
\]
From the Lax-Milgram theorem and the standard interior elliptic estimates, $T_{\mu}(u)$ is well defined and belongs to $X\cap C^{2,\alpha}(\R)$ for every $u\in X$. The operator $T_{\mu}$ is continuous, and compact in the sense that it maps bounded subsets of $X$ into subsets of $X$ 
with compact closure. Furthermore, again from the elliptic estimates, the map 
$\mu\mapsto T_{\mu}$ is continuous  locally uniformly in $X$, in the sense that, 
for every $\mu>0$, every sequence $(\mu_n)_{n\in\N}$ converging to $\mu$ and 
every bounded subset $B$ of $X$, one has 
\[
\hbox{$\sup_{u\in B}\|T_{\mu_n}(u)-T_{\mu}(u)\|_X\to0$ as $n\to+\infty$.}
\]

Notice also that the constant function $1\in X$ is a fixed point of $T_{\mu}$, that is, $T_{\mu}(1)=1$, for every~$\mu>0$. In order to get the conclusion of Theorem~\ref{th1}, our goal 
is to show that, when~$\mu>0$ is large enough, there are non-constant positive fixed points of $T_{\mu}$ in $X$, that is, positive non-constant solutions $u\in X$ of 
\[
u-T_{\mu}(u)=0.
\]
We will do so by evaluating the Leray-Schauder topological degree of the map $I-T_{\mu}$ at the point~$0$ in some suitable open subsets of $X$ and by using a homotopy argument from small values of the parameter $\mu$ to large values of~$\mu$. Here,~$I$ denotes the identity map.

The first step is to show that the constant function $u\equiv 1$ is the only positive fixed point of~$T_{\mu}$ in~$X$ when $\mu>0$ is small enough.
\begin{lem}\label{lem2}
There is $\underline{\mu}>0$ such that, for every $\mu\in(0,\underline{\mu})$, the constant function~$u\equiv 1$ is the only positive fixed point of~$T_{\mu}$ in~$X$, that is the only positive solution of~$(\ref{stat})$ in $X$. 
\end{lem}

\noindent{\bf{Proof.}} In Theorem~1.1 of~\cite{bnpr}, the stronger conclusion that the 
constant function~$1$ is the only positive bounded solution of~(\ref{stat}) for small 
$\mu>0$ was proved, under some slightly stronger smoothness assumptions on $\phi$. 
In the present Lemma~\ref{lem2}, uniqueness holds with slightly different assumptions on 
$\phi$, but is shown here only in the space~$X$, so the proof is shorter. 
It is presented for the sake of completeness.

Assume that the conclusion is not true. Then there is a sequence $(\mu_n)_{n\in\N}$ of positive 
real numbers such that $\mu_n\to0$ as $n\to+\infty$, and a sequence $(u_n)_{n\in\N}$ in $X$ 
of solutions of~(\ref{stat}) with $\mu=\mu_n$, such that $0<u_n$ in $\R$ and $u_n\not\equiv 1$. 
From Lemma~\ref{lem1}, there is a constant $M$ such that
\[
0<u_n\le M\hbox{ in }\R,
\]
for every $n\in\N$. If $\max_{\R}u_n\le1$, then $\phi*u_n\le1$ and $u_n''\le0$ in $\R$, 
implying that $u_n$ is a positive constant, which would then have to be equal to $1$ by~(\ref{stat}). Therefore,
\be\label{maxun}
\max_{\R}u_n>1\hbox{ for every }n\in\N.
\ee\par
As done in the proof of Lemma~\ref{lem1}, it follows from~(\ref{stat}) that the 
$L$-periodic functions $u_n$ are bounded in $C^k(\R)$ for every $k\in\N$, and therefore converge, 
as $n\to+\infty$, up to extraction of a subsequence,  to a $C^{\infty}(\R)$ solution $u_{\infty}$ of 
\[
u_{\infty}''=0,
\]
such that $0\le u_{\infty}\le M$ in~$\R$. As a consequence, the function~$u_{\infty}$ is constant and this constant is such that $1\le u_{\infty}\le M$ due to~(\ref{maxun}). If $u_{\infty}>1$, then,
as all $u_n$ are $L$-periodic, we have $u_n>1$ in~$\R$ for~$n$ large enough, whence 
$\phi*u_n>1$ and $u''_n>0$ in~$\R$ for~$n$ large enough, which is impossible since all $u_n$ are bounded. We conclude that
\[
u_{\infty}\equiv 1.
\]
Next, write $u_n=1+v_n$ and $w_n=v_n/\|v_n\|_{L^{\infty}(\R)}$. Since $\max_{\R}u_n>1$, one has $\max_{\R}v_n>0$ and $w_n$ is well defined, for every $n\in\N$. Furthermore, $v_n\to0$ as $n\to+\infty$ in $C^k(\R)$ for every $k\in\N$. The functions $w_n$ are $L$-periodic, with $\max_{\R}w_n>0$ and $\|w_n\|_{L^{\infty}(\R)}=1$. They obey
\[
w_n''-\mu_n\,u_n\,(\phi*w_n)=0\hbox{ in }\R.
\]
From the standard elliptic estimates, we know that the functions $w_n$ converge as $n\to+\infty$, up to extraction of a subsequence, in $C^k(\R)$ for every $k\in\N$, to a $C^{\infty}(\R)$ solution $w_{\infty}$ of 
\[
w_{\infty}''=0,
\]
such that $\max_{\R}w_{\infty}\ge0$ and $\|w_{\infty}\|_{L^{\infty}(\R)}=1$. It follows that 
$w_{\infty}\equiv 1$. For $n$ large enough, one gets that $w_n>0$ in $\R$, once again 
since $w_n$ are $L$-periodic, whence
\[
u_n=1+\|v_n\|_{L^{\infty}(\R)}w_n>1\hbox{ in }\R.
\]
Finally, it follows that $\phi*u_n>1$ and $u''_n>0$ in $\R$ for $n$ large enough, which is impossible since $u_n$ is bounded. One has then reached a contradiction and the proof 
of Lemma~\ref{lem2} is thereby complete.~\hfill$\Box$

\subsubsection*{Stability of the fixed points}

Let us now study the stability of the fixed point $1$ of $T_{\mu}$ when $\mu>0$ is small or large. It follows from elementary calculations and standard elliptic estimates that $T_{\mu}$ is Fr\'echet-differentiable everywhere in $X$ and that the Fr\'echet-derivative $DT_{\mu}(1)$ at the point $1$ is given as follows: for every $u\in X$, $w=DT_{\mu}(1)(u)$ is the unique solution of
$$-w''+w=u-\mu\,\phi*u\hbox{ in }\R.$$
Notice that the linear operator $DT_{\mu}(1):X\to X$ is compact.\par
The next two lemmas are concerned with the stability of the fixed point $1$ of~$T_{\mu}$ when~$\mu>0$ is small and when $\mu>0$ is large.

\begin{lem}\label{lem3} $($Small $\mu)$
There is $\underline{\underline{\mu}}>0$ such that, for every $\mu\in(0,\underline{\underline{\mu}})$, all eigenvalues $\lambda$ of $DT_{\mu}(1)$ satisfy
$|\lambda|<1$.
\end{lem}
{\bf Proof.} Let $\underline{\underline{\mu}}$ be   given by
\[
\underline{\underline{\mu}}=\min\Big(\frac{4\pi^2}{L^2},\frac{1}{2}\Big),
\]
and let us check that the conclusion holds with this value $\underline{\underline{\mu}}$. To this 
end, let $\mu\in(0,\underline{\underline{\mu}})$ and let $u$ be an eigenvector of $DT_{\mu}(1)$ 
with the eigenvalue $\lambda$. The function $u\in X$ solves
\be\label{lambdau}
-(\lambda u)''+\lambda u=u-\mu\,\phi*u\hbox{ in }\R.
\ee
Let $(a_k)_{k\in\N}$ be the Fourier coefficients of $u$, defined by
\be\label{defam}
a_k=\frac{2}{L}\int_0^Lu(x)\,\cos\Big(\frac{2\pi kx}{L}\Big)\,dx,\ \ k\in\N.
\ee
The Fourier coefficients $b_k$ of $\phi*u$ are given by
\[
b_k=\frac{2}{L}\int_0^L(\phi*u)(x)\,
\cos\Big(\frac{2\pi kx}{L}\Big)\,dx=
\hat{\phi}\Big(\frac{k}{L}\Big)\,a_k,\ \ k\in\N,
\]
since $\phi$ is even and real-valued.

Since $u$ is even, $L$-periodic and does not vanish identically, there is $m_0\in\N$ such 
that $a_{m_0}\neq 0$. If $\lambda$ were equal to $0$, then, in particular, the $m_0$-th Fourier coefficient of the right-hand side of~(\ref{lambdau}) would vanish, that is
\[
a_{m_0}-\mu\,\hat{\phi}\Big(\frac{m_0}{L}\Big)\,a_{m_0}=0,
\]
whence 
\[
1-\mu\,\hat{\phi}(m_0/L)=0,
\]
as $a_{m_0}\neq0$. Thus, we would have 
\[
1=\mu\,\hat{\phi}(m_0/L)\le\underline{\underline{\mu}}\|\phi\|_{L^1(\R)}=
\underline{\underline{\mu}}\le1/2,
\]
which is impossible. As a consequence, $\lambda\neq 0$. It follows then from~(\ref{lambdau}) 
that the function $u$ is actually of class $C^{\infty}(\R)$ and that
\[
\Big(\frac{4\pi^2k^2}{L^2}+1\Big)\lambda\,a_k=a_k-\mu\,\hat{\phi}
\Big(\frac{k}{L}\Big)a_k\ \hbox{ for all }k\in\N.
\]
In particular, since $a_{m_0}\neq 0$, one gets that
\[
\lambda\Big(\frac{4\pi^2m_0^2}{L^2}+1\Big)=1-\mu\,\hat{\phi}\Big(\frac{m_0}{L}\Big).
\]
If $m_0=0$, then $\lambda=1-\mu$, whence $\lambda\in(1-\underline{\underline{\mu}},1)\subset(1/2,1)$. If $m_0\ge 1$, then
\[
|\lambda|\,\Big(\frac{4\pi^2m_0^2}{L^2}+1\Big)=\Big|1-\mu\,
\hat{\phi}\Big(\frac{m_0}{L}\Big)\Big|<1+\underline{\underline{\mu}}\le1+\frac{4\pi^2}{L^2},
\]
whence $|\lambda|<1$. The proof of Lemma~\ref{lem3} is thereby complete.\hfill$\Box$

\begin{lem}\label{lem4} $($Large $\mu)$
There is $\mu^*>0$ such that, for every $\mu>\mu^*$, the operator $DT_{\mu}(1)$ has 
an eigenvalue~$\lambda_{\mu}>1$, with ${\rm{dim}}\ker(DT_{\mu}(1)-\lambda_{\mu}I)^n=1$ 
for all $n\ge 1$. Moreover, all other eigenvalues~$\lambda$ of~$DT_{\mu}(1)$  satisfy $\lambda<1$.
\end{lem}

\noindent{\bf{Proof.}} Recall that we assume that there exists $k_0\in\N$ such that $\hat{\phi}(k_0/L)<0$ 
and $\hat{\phi}(k/L)\ge0$ for all $k\in\N\backslash\{k_0\}$. As $\hat\phi(0)=1$, we have $k_0\ge 1$.
Let $\mu_*>0$ be defined by
\begin{equation}\label{defmu*}
\mu^*=\frac{4\pi^2k_0^2}{L^2}\times
\frac{1}{\Big|\hat{\phi}\big(\displaystyle\frac{k_0}{L}\big)\Big|}>0
\end{equation}
and let us show that the conclusion of Lemma~\ref{lem4} holds with this choice of $\mu^*$.
Let $\mu$ be any real number such that $\mu>\mu^*$, and set
\begin{equation}\label{deflambdamu}
\lambda_{\mu}=\frac{1-\mu\,\hat{\phi}\big(\displaystyle
\frac{k_0}{L}\big)}{\displaystyle\frac{4\pi^2k_0^2}{L^2}+1}>1,
\end{equation}
by~(\ref{defmu*}). The function $u_{\mu}\in X$ given by
\[
u_{\mu}(x)=\cos\Big(\frac{2\pi k_0x}{L}\Big)\ \hbox{ for every }x\in\R,
\]
solves
\[
-\lambda_{\mu}u_{\mu}''+\lambda_{\mu}u_{\mu}=
\Big(1-\mu\,\hat{\phi}\big(\frac{k_0}{L}\big)\Big)u_{\mu}=u_{\mu}-\mu\,\phi*u_{\mu}\hbox{ in }\R.
\]
In other words, $DT_{\mu}(1)(u_{\mu})=\lambda_{\mu}u_{\mu}$, that is $u_{\mu}$ is an eigenvector of $DT_{\mu}(1)$ with the eigenvalue~$\lambda_{\mu}>1$.

Moreover, if  $u\in X$ satisfies $DT_{\mu}(1)(u)=\lambda_{\mu}u$, that is
\begin{equation}\label{lambdamuu}
-\lambda_{\mu}u''+\lambda_{\mu}u=u-\mu\,\phi*u\hbox{ in }\R,
\end{equation}
then, since $\lambda_{\mu}\neq0$, the function $u$ is actually of class $C^{\infty}(\R)$ and its Fourier coefficients $a_k$ given by~(\ref{defam}) satisfy
\be\label{aku}
\Big(\frac{4\pi^2k^2}{L^2}+1\Big)\lambda_{\mu}\,a_k=\Big(1-\mu\,\hat{\phi}\big(\frac{k}{L}\big)\Big)\,a_k
\ee
for every $k\in\N$. For every $k\in\N\backslash\{k_0\}$, there holds
$$\Big(\frac{4\pi^2k^2}{L^2}+1\Big)\lambda_{\mu}\ge\lambda_{\mu}>1\ge 1-\mu\,\hat{\phi}\big(\frac{k}{L}\big)$$
by the assumption~(\ref{hypphi}). Therefore, $a_k=0$ for every~$k\in\N\backslash\{k_0\}$. 
Since $u$ is even, it is then a multiple of the function $u_{\mu}(x)=\cos(2\pi k_0x/L)$.

Next, let $u\in X$ be such that 
\[
(DT_{\mu}(1)-\lambda_{\mu}I)^2(u)=0,
\]
that is $v=DT_{\mu}(1)u-\lambda_{\mu}u$ is in the eigenspace of $DT_{\mu}(1)$ with 
eigenvalue $\lambda_{\mu}$. From the previous paragraph, one knows that 
$v=\gamma u_{\mu}$ for some $\gamma\in\R$. In other words, 
\[
DT_{\mu}(1)(u)=\lambda_{\mu}u+\gamma u_{\mu}.
\] 
In particular, the function $u$ is of class~$C^{\infty}(\R)$ and solves
\[
-(\lambda_{\mu}u+\gamma u_{\mu})''+(\lambda_{\mu}u+\gamma u_{\mu})=u-\mu\,\phi*u.
\]
Therefore, the Fourier coefficients $a_k$ of $u$, given by~(\ref{defam}), satisfy~(\ref{aku}) for every $k\in\N\backslash\{k_0\}$. As in the previous paragraph, it follows that $a_k=0$ for every $k\in\N\backslash\{k_0\}$, hence, $u$ is a multiple of the function $u_{\mu}$.

As a consequence, ${\rm{dim}}\,\ker(DT_{\mu}(1)-\lambda_{\mu}I)^n=1$ for every $n\ge 1$, and 
the algebraic and geometric multiplicity of the eigenvalue $\lambda_{\mu}$ of $DT_{\mu}(1)$ 
is equal to $1$.

It only remains to show that the other eigenvalues of $DT_{\mu}(1)$ are all less than $1$. 
To do so, let~$\lambda$ be any real number such that $\lambda\neq\lambda_{\mu}$ and 
$\lambda\ge 1$, and let $u\in X$ be such that $DT_{\mu}(1)(u)=\lambda u$. One has to prove 
that $u$ is necessarily equal to $0$. The function $u$ is of class $C^{\infty}(\R)$, 
it solves~(\ref{lambdamuu}) with the parameter $\lambda$ and its Fourier coefficients $a_k$ 
given by~(\ref{defam}) satisfy~(\ref{aku}) with the parameter~$\lambda$, 
that is,
\[
\Big(\frac{4\pi^2k^2}{L^2}+1\Big)\lambda\,a_k=
\Big(1-\mu\,\hat{\phi}\big(\frac{k}{L}\big)\Big)\,a_k\ \hbox{ for every }k\in\N.
\]
Given the definition~(\ref{deflambdamu}) of~$\lambda_{\mu}$, and given 
that $\lambda\neq\lambda_{\mu}$, it follows that $a_{k_0}=0$. Furthermore, 
for~$k=0$, one has $\lambda\ge 1>1-\mu=1-\mu\hat{\phi}(0)$, whence $a_0=0$. 
Lastly, for every $k\in\N\backslash\{0,k_0\}$, there holds
\[
\Big(\frac{4\pi^2k^2}{L^2}+1\Big)\lambda>\lambda\ge 1\ge 1-\mu\hat{\phi}\big(\frac{k}{L}\big),
\]
whence $a_k=0$. As a consequence, $a_k=0$ for every $k\in\N$, that is $\lambda$ cannot 
be an eigenvalue of~$DT_{\mu}(1)$. That completes the proof of Lemma~\ref{lem4}.\hfill$\Box$

\begin{rem}{\rm The fact that $X$ only contains even $L$-periodic functions forces the operator $DT_{\mu}(1)$, for large $\mu$, to have only one eigendirection associated to an eigenvalue larger than $1$. That will yield the explicit value, namely~$-1$, of the index of the map $I-T_{\mu}$ at $0$ in a small neighborhood around the point $u=1$ and will finally lead to the existence of other $($than $u=1)$ fixed points of~$T_{\mu}$ in~$X$ for large~$\mu$.}
\end{rem}

\subsubsection*{End of the proof of Theorem~\ref{th1}} 

Let $\underline{\mu}>0$ and $\underline{\underline{\mu}}>0$ be as in 
Lemmas~\ref{lem2} and~\ref{lem3}, and set
\[
\mu_*=\min\Big(\underline{\mu},\underline{\underline{\mu}}\Big)>0.
\]
Let $\mu^*>0$ be as in Lemma~\ref{lem4}. From Lemmas~\ref{lem3} and~\ref{lem4}, it follows that $\mu_*\le\mu^*$.

Let now $\mu_0$ and $\tilde{\mu}$ be any two real numbers such that
\[
0<\mu_0<\mu_*\le\mu^*<\tilde{\mu}.
\]
From Lemma~\ref{lem1}, there are two positive constants $0<m\le M$ such that every positive 
bounded solution~$u$ of~(\ref{stat}), with parameter $\mu\in[\mu_0,\tilde{\mu}]$, 
satisfies
\[
0<m\le u(x)\le M\hbox{ for all }x\in\R.
\]
Notice that $u\equiv 1$ is a solution of~(\ref{stat}) for every $\mu$, whence $0<m\le 1\le M$.

Let $\Omega$ be the non-empty open bounded subset of $X$ defined by
\[
\Omega=\Big\{u\in X,\ \frac{m}{2}<\min_{\R}u\le\max_{\R}u<2M\Big\}.
\]
It follows from Lemma~\ref{lem1} and the choice of $m$ and $M$ that, for every $\mu\in[\mu_0,\tilde{\mu}]$ and for every~$u\in\partial\Omega$, there holds
\[
u-T_{\mu}u\neq0,
\]
for 
otherwise $u$ would be a positive solution of~(\ref{stat}) with the value $\mu$, whence 
$m\le u\le M$ in $\R$ and $u\not\in\partial\Omega$. Therefore, for every $\mu\in[\mu_0,\tilde{\mu}]$, the Leray-Schauder topological degree 
$\hbox{deg}(I-T_{\mu},\Omega,0)$ of the map $I-T_{\mu}$ in the set $\overline{\Omega}$ at the point $0$ is well defined. Furthermore, by the continuity of $T_{\mu}$ with respect to $\mu$ 
in the local uniform sense in $X$, this degree   does not depend on~$\mu\in[\mu_0,\tilde{\mu}]$.

Let us now compute the degree. Lemmas~\ref{lem2} and~\ref{lem3} imply that $\hbox{deg}(I-T_{\mu_0},\Omega,0)$ is equal to the index $\hbox{ind}(I-T_{\mu_0},1,0)$ of $I-T_{\mu_0}$ at the point $u=1$, namely $\hbox{ind}(I-T_{\mu_0},1,0)=1$. In other words, $\hbox{deg}(I-T_{\mu_0},\Omega,0)=1$, whence
\begin{equation}\label{degree1}
\hbox{deg}(I-T_{\tilde{\mu}},\Omega,0)=1.
\end{equation}

On the other hand, it follows from Lemma~\ref{lem4} that the operator $I-DT_{\tilde{\mu}}(1)$ is 
one-to-one. Since $I-DT_{\tilde{\mu}}(1)$ is a Fredholm operator of index $0$, it is 
invertible and the point $u=1$ is an isolated zero of $I-T_{\tilde{\mu}}$ in $X$. Finally, $\hbox{ind}(I-T_{\tilde{\mu}},1,0)=-1$ since $I-DT_{\tilde{\mu}}(1)$ has only one direction associated to a negative eigenvalue. In other words,
\begin{equation}\label{degree2}
\hbox{deg}(I-T_{\tilde{\mu}},B_{\epsilon}(1),0)=\hbox{ind}(I-T_{\tilde{\mu}},1,0)=-1
\end{equation}
for $\epsilon>0$ small enough, where $B_{\epsilon}(1)$ denotes the open ball of center $1$ and radius $\epsilon$ in $X$. But
\[
\hbox{deg}(I-T_{\tilde{\mu}},\Omega,0)=\hbox{deg}(I-T_{\tilde{\mu}},B_{\epsilon}(1),0)+\hbox{deg}(I-T_{\tilde{\mu}},\Omega\backslash B_{\epsilon}(1),0),
\]
for $\epsilon>0$ small enough. Together with~(\ref{degree1}) and~(\ref{degree2}), one finally concludes that, for $\epsilon>0$ small enough,
\[
\hbox{deg}(I-T_{\tilde{\mu}},\Omega\backslash B_{\epsilon}(1),0)=2.
\]
For $\epsilon>0$ small enough, there exists then a solution $u$ of $u-T_{\tilde{\mu}}(u)=0$ in $\Omega\backslash B_{\epsilon}(1)$, that is a non-constant positive bounded and $L$-periodic solution~$u$ of~(\ref{stat}) with the parameter $\tilde{\mu}$ (remember that $1$ is the only positive constant function solving~(\ref{stat})). The proof of Theorem~\ref{th1} is thereby complete.\hfill$\Box$


\SE{Bounds for the solutions of the Cauchy problem}\label{sec3}

Here, we prove Theorem~\ref{th2}. We consider the Cauchy
problem~(\ref{cauchy}) with a parameter $\mu>0$ and a nonnegative
initial condition $u_0\in L^{\infty}(\R)$. We point out that we do not
assume that $u_0$ is periodic.  The maximum principle and standard
parabolic estimates imply that the solution $u$ exists for all times
$t\in(0,+\infty)$, is classical in~$(0,+\infty)\times\R$ (and is even
of class $C^{\infty}$ in this set) and satisfies
\be\label{localbounds} 0\le u(t,x)\le 
e^{\mu t}\,\|u_0\|_{L^{\infty}(\R)}\ \hbox{ for every }t>0\hbox{ and
}x\in\R.  
\ee
Our task is to improve (\ref{localbounds}) to a uniform in time
estimate. 

Let $\sigma>0$ be as in~(\ref{hypphi0}) and define the local average
on the scale $\sigma$:
\[
v(t,x)=\int_{x-\sigma/2}^{x+\sigma/2}u(t,y)\,dy\ \hbox{ for
}(t,x)\in[0,+\infty)\times\R.
\]
The function $v$ is of class $C^{\infty}((0,+\infty)\times\R)$,
continuous in 
$[0,+\infty)\times\R$ and satisfies a version of the upper bound (\ref{localbounds}):
$$0\le v(t,x)\le \sigma\,e^{\mu t}\|u_0\|_{L^{\infty}(\R)}\ \hbox{ for every }(t,x)\in[0,+\infty)\times\R.$$
Furthermore, the function $v$ obeys
\[
v_t(t,x)-v_{xx}(t,x)=\mu\int_{x-\sigma/2}^{x+\sigma/2}u(t,y)\,(1-(\phi*u)(t,y))\,dy,
\]
for every $(t,x)\in(0,+\infty)\times\R$, and, since the right-hand
side of  the above equation belongs to~$L^{\infty}((a,b)\times\R)$ for
every $0\le a<b<+\infty$, the function
$t\mapsto\|v(t,\cdot)\|_{L^{\infty}(\R)}$ is actually continuous on
$[0,+\infty)$. 

Owing to~(\ref{hypphi0}), let now $\eta>0$ be such that
\be\label{defeta} \phi\ge\eta>0\hbox{ a.e. in }(-\sigma,\sigma), \ee
and let $M$ be any positive real number such that 
\begin{equation}\label{defM}
M>\max\Big(\sigma\,\|u_0\|_{L^{\infty}(\R)},\frac{1}{\eta}\Big).  
\end{equation}
We will show that $\|v(t,\cdot)\|_{L^{\infty}(\R)}\le M$ for all
$t>0$, by contradiction. Assume that this is false. Since
$t\mapsto\|v(t,\cdot)\|_{L^{\infty}(\R)}$ is continuous on
$[0,+\infty)$, and
\[
\|v(0,\cdot)\|_{L^{\infty}(\R)}\le\sigma\|u_0\|_{L^{\infty}(\R)}<M,
\]
there exists $t_0>0$ such that
$\|v(t_0,\cdot)\|_{L^{\infty}(\R)}=M$ and
$\|v(t,\cdot)\|_{L^{\infty}(\R)}<M$ for all $t\in[0,t_0)$. Since $v$
is nonnegative, there exists then a sequence of real numbers
$(x_n)_{n\in\N}$ such that $v(t_0,x_n)\to M$ as $n\to+\infty$. As
usual, we define the translates
\[
u_n(t,x)=u(t,x+x_n)\ \hbox{ and }\ v_n(t,x)=v(t,x+x_n),
\]
for $n\in\N$ and $(t,x)\in(0,+\infty)\times\R$. From standard
parabolic estimates, the sequences~$(u_n)_{n\in\N}$
and~$(v_n)_{n\in\N}$ are bounded in $C^k_{loc}((0,+\infty)\times\R)$
for every $k\in\N$ and converge in these spaces, up to extraction of a subsequence, to
some nonnegative functions $u_{\infty}$ and $v_{\infty}$ of class
$C^{\infty}((0,+\infty)\times\R)$, such that
\[
v_{\infty}(t,x)=\int_{x-\sigma/2}^{x+\sigma/2}u_{\infty}(t,y)\,dy,
\]
and
\[
(v_{\infty})_t(t,x)=(v_{\infty})_{xx}(t,x)+\mu\int_{x-\sigma/2}^{x+\sigma/2}u_{\infty}(t,y)\,
(1-(\phi*u_{\infty})(t,y))\,dy,
\]
for every $(t,x)\in(0,+\infty)\times\R$. The passage to the limit in
the integral term is possible due to the local uniform convergence of
$u_n$ to $u_{\infty}$ in $(0,+\infty)\times\R$ and to the
local-in-time bounds~(\ref{localbounds}). Furthermore, we have
\[
0\le
v_{\infty}(t,x)\le M,
\]
for every $0<t\le t_0$ and $x\in\R$,
and~$v_{\infty}(t_0,0)=M$. Therefore, we have
\[
\hbox{$(v_{\infty})_t(t_0,0)\ge0$ and
$(v_{\infty})_{xx}(t_0,0)\le 0$,}
\]
whence,
\[
\int_{-\sigma/2}^{\sigma/2}u_{\infty}(t_0,y)\,(1-(\phi*u_{\infty})(t_0,y))\,dy\ge0.
\]
If 
\begin{equation}\label{phistar}
\hbox{$(\phi*u_{\infty})(t_0,\cdot)>1$ everywhere in
$[-\sigma/2,\sigma/2]$},
\end{equation}
then the continuous function
\[
U=u_{\infty}(t_0,\cdot)\,(1-(\phi*u_{\infty})(t_0,\cdot))
\]
would be
nonpositive on $[-\sigma/2,\sigma/2]$. Since its integral over
$[-\sigma/2,\sigma/2]$ is nonnegative, the function~$U$ would be
identically equal to zero on $[-\sigma/2,\sigma/2]$. Moreover, it would
then follow from (\ref{phistar}) that $u_{\infty}(t_0,\cdot)=0$ on $[-\sigma/2,\sigma/2]$,
whence $v_{\infty}(t_0,0)=0$, contradicting the assumption 
that~$v_{\infty}(t_0,0)=M>0$. Therefore, there is a real number
$y_0\in[-\sigma/2,\sigma/2]$ such
that
\[
(\phi*u_{\infty})(t_0,y_0)\le1.
\]
Since both functions $\phi$ and
$u_{\infty}$ are nonnegative, one gets from~(\ref{defeta}) that
$$\baa{rcl}
1\ge(\phi*u_{\infty})(t_0,y_0) & \!\!\ge\!\! & \displaystyle\int_{-\sigma}^{\sigma}\phi(y)\,u_{\infty}(t_0,y_0-y)\,dy\vspace{3pt}\\
& \!\!\ge\!\! & \displaystyle\eta\int_{-\sigma}^{\sigma}u_{\infty}(t_0,y_0-y)\,dy\ge\eta\int_{-\sigma/2}^{\sigma/2}u_{\infty}(t_0,y)\,dy=\eta\,v_{\infty}(t_0,0)=\eta\,M.\eaa$$
This contradicts the definition~(\ref{defM}) of the constant $M$.

We conclude that $\|v(t,\cdot)\|_{L^{\infty}(\R)}\le M$ for all $t\ge
0$. Since $u$ is nonnegative, this means that 
\begin{equation}\label{boundsv}
0\le\int_{x-\sigma/2}^{x+\sigma/2}u(t,y)\,dy\le M 
\end{equation}
for every
$t\ge0$ and $x\in\R$. To get a global bound for $u$ itself, we just fix an
arbitrary time $s\ge1$ and infer from the maximum principle that
\[
0\le u(s,x)\le w(s,x)\ \hbox{ for every }x\in\R.
\]
Here, $w$ is the solution of the equation
\[
w_t=w_{xx}+\mu\,w
\]
with the initial condition at time $s-1$ given by
$w(s-1,\cdot)=u(s-1,\cdot)$. It follows then
from~(\ref{boundsv}) that, for every $x\in\R$,
$$0\le u(s,x)\le e^{\mu}\int_{-\infty}^{+\infty}\frac{e^{-y^2/4}}{\sqrt{4\pi}}\,u(s-1,x-y)\,dy\le\frac{2\,M\,e^{\mu}}{\sqrt{4\pi}}\sum_{k\in\N}e^{-k^2\sigma^2/4}<+\infty.$$
Together with the local-in-time bounds~(\ref{localbounds}), this
implies that $u$ is globally bounded. The proof of Theorem~\ref{th2}
is complete.\hfill$\Box$


\SE{Spreading speed for the Cauchy problem~(\ref{cauchy})}\label{sec4}

This section is devoted to the proof of Theorem~\ref{th3} on the
asymptotic spreading speed for the solutions of the Cauchy
problem~(\ref{cauchy}). The proof of the upper bound~(\ref{spreadsup})
is immediate simply by comparing the solution $u$ of (\ref{cauchy}) to that of the (local)
linear heat equation 
\begin{equation}\label{eqlin} 
v_t=v_{xx}+\mu\,v.  
\end{equation} 
The proof of the lower bound~(\ref{spreadinf}) is more involved and,
when compared to the proofs given in~\cite{aw,bhn1,bhn2} for the
analogous local problems, it requires additional arguments due to the
nonlocality of~(\ref{eq}), and the lack of comparison principle for the
nonlocal equation~(\ref{eq}). However, the maximum principle can be
applied when $\phi*u$ is small and the equation~(\ref{eq}) will then
be compared in some suitable bounded boxes to a local linear equation
close to~(\ref{eqlin}).

\subsubsection*{{Proof of the upper bound }} 

We assume here that $u_0$ is compactly supported. Hence, there is
$R>0$ such that $u_0(x)=0$ for a.e. $|x|\ge R$. Since $u(t,x)\ge0$ for
all $t>0$ and $x\in\R$, one has
\[
\mu\,u(t,x)\,(1-(\phi*u)(t,x))\le\mu\,u(t,x)
\]
for all $t>0$ and $x\in\R$. Let $v$ denote the solution
of~(\ref{eqlin}) for $t>0$ with initial condition $u_0$ at~$t=0$. It
follows from the maximum principle that $0\le u(t,x)\le v(t,x)$ for
all $t>0$ and $x\in\R$, whence 
\begin{equation}\label{ineqsup} 0\le
u(t,x)\le\frac{e^{\mu t}}{\sqrt{4\pi
    t}}\int_{-R}^Re^{-(x-y)^2/(4t)}\,u_0(y)\,dy.  
\end{equation} 
Let $c$ be any arbitrary real number such that $c>2\sqrt{\mu}$. For
all $t\ge R/c$ and for all $|x|\ge ct$, one has
\[
0\le u(t,x)\le\frac{e^{\mu t}\|u_0\|_{L^{\infty}(\R)}}{\sqrt{4\pi
    t}}\int_{-R}^Re^{-(ct-R)^2/(4t)}dy,
\]
which immediately yields~(\ref{spreadsup}) (notice that, for every
$t>0$, the maximum of $u(t,\cdot)$ on the
set~$(-\infty,-ct]\cup[ct,+\infty)$ is reached since $u(t,\cdot)$ is
continuous, nonnegative, and converges to~$0$ at~$\pm\infty$
from~(\ref{ineqsup})).

\subsubsection*{{Proof of the lower bound}} 

Assume for the sake of contradiction that the conclusion~(\ref{spreadinf}) does not
hold. Then, since $u$ is nonnegative, there is $c$ such that
\[
0\le c<2\sqrt{\mu}
\]
and there are two sequences $(t_n)_{n\in\N}$ in $(0,+\infty)$ and $(x_n)_{n\in\N}$ in~$\R$ such that
\begin{equation}\label{utnxn}\left\{\baa{l}
|x_n|\le c\,t_n\hbox{ for all }n\in\N,\vspace{3pt}\\
t_n\to+\infty\hbox{ and }u(t_n,x_n)\to0\hbox{ as }n\to+\infty.\eaa\right.
\end{equation}
We set
\begin{equation}\label{defcn}
c_n=\frac{x_n}{t_n}\,\in[-c,c].
\end{equation}
Up to extraction of a subsequence, one can assume that $c_n\to
c_{\infty}\,\in[-c,c]$ as $n\to+\infty$.

Furthermore, for every $n\in\N$ and $(t,x)\in(-t_n,+\infty)\times\R$,
we define the shifted functions
\[
u_n(t,x)=u(t+t_n,x+x_n).
\]
From Theorem~\ref{th2}, the sequence
$(\|u_n\|_{L^{\infty}((-t_n,+\infty)\times\R)})_{n\in\N}$ is
bounded. Therefore, the standard parabolic estimates imply that the
functions $u_n$ converge in $C^{1,2}_{loc}(\R\times\R)$, up to
extraction of a subsequence, to a classical bounded solution
$u_{\infty}$ of
\[
(u_{\infty})_t=(u_{\infty})_{xx}+\mu\,u_{\infty}\,(1-\phi*u_{\infty})\
\hbox{ in }\R\times\R,
\]
such that $u_{\infty}\ge0$ in $\R\times\R$ and $u_{\infty}(0,0)=0$. By
viewing $\mu\,(1-\phi*u_{\infty})$ as a coefficient
in~$L^{\infty}(\R\times\R)$, it follows from the strong parabolic
maximum principle and the uniqueness of the solutions of the Cauchy
problem that $u_{\infty}(t,x)=0$ for all~$(t,x)\in\R\times\R$ (the
limit $u_{\infty}$ being unique, one can then infer that the whole
sequence $(u_n)_{n\in\N}$ converges to $0$ in
$C^{1,2}_{loc}(\R\times\R)$). As a consequence, the nonnegative
functions $v_n$ defined in $(-t_n,+\infty)\times\R$ by
\[
v_n(t,x)=u_n(t,x+c_nt)=u(t+t_n,x+c_n(t+t_n)) 
\]
converge to $0$ locally uniformly in $\R\times\R$, due to the
boundedness of the speeds $c_n$ defined in~(\ref{defcn}). Hence the
nonnegative functions $\phi*v_n$ also converge to $0$ locally
uniformly in $\R\times\R$, since the
sequence~$(\|v_n\|_{L^{\infty}((-t_n,+\infty)\times\R)})_{n\in\N}$ is
bounded.

Let us now fix some parameters which are independent of $n$. First,
let $\delta>0$ be such that
\begin{equation}\label{defdelta} 
\mu\,(1-\delta)\ge\frac{c^2}{4}+\delta ,
\end{equation}
and also let
$R>0$ be such that 
\begin{equation}
\label{defR} \frac{\pi^2}{4R^2}\le\delta.  
\end{equation}
Since $u(1,\cdot)$ is continuous from parabolic regularity, and is
positive in $\R$ from the strong parabolic maximum principle, there is
$\eta>0$ such that \be\label{defeta2} u(1,x)\ge\eta>0\ \hbox{ for all
}|x|\le R+c.  \ee 

Without loss of generality, one can assume that $t_n>1$ for every
$n\in\N$.  Since $\phi*v_n\to0$ locally
uniformly in $\R\times\R$ as $n\to+\infty$, for $n$ sufficiently large ($n\ge N$) 
we may define 
\[
t^*_n=\inf\Big\{t\in[-t_n+1,0];\ (0\le)\,\phi*v_n\le\delta\hbox{ in
}[t,0]\times[-R,R]\Big\},~~n\ge N,
\]
with $\delta$ and $R$ as in~(\ref{defdelta}) and~(\ref{defR}), and we may assume that $t^*_n<0$. Furthermore, for every $n\ge N$, by continuity of $\phi*v_n$ in
$(-t_n,+\infty)\times\R$, the infimum is a minimum in the definition
of $t^*_n$ and 
\begin{equation}\label{delta2} 0\le\phi*v_n\le\delta\hbox{ in
}[t^*_n,0]\times[-R,R].  
\end{equation} 
On the other hand, we have
\[
v_n(-t_n+1,x)=u(1,x+c_n)\ge\eta\hbox{ for all }|x|\le R,
\]
from~(\ref{defcn}) and~(\ref{defeta2}), for all $n\in\N$. Therefore,
again by continuity of $\phi*v_n$ in $(-t_n,+\infty)\times\R$ and by
minimality of $t^*_n$, for each $n\ge N$ there is the following dichotomy:
\be\label{dichotomy}\left\{\baa{ll}
\hbox{either} & \Big(t^*_n\!>\!-t_n\!+\!1\hbox{ and }
\displaystyle{\mathop{\max}_{[-R,R]}}(\phi*v_n)(t^*_n,\cdot)=\delta\Big)\vspace{3pt}\\
  \hbox{or} & \Big(t^*_n\!=\!-t_n\!+\!1\hbox{ and
  }\displaystyle{\mathop{\min}_{[-R,R]}}v_n(t^*_n,\cdot)\ge\eta\Big).\eaa\right.  
\ee

Next, we claim that there exists $\rho>0$ such that
\begin{equation}\label{claimrho} 
\min_{[-R,R]}v_n(t^*_n,\cdot)\ge\rho>0\ \hbox{ for
  all }n\ge N.
\end{equation} 
Note that this claim would be immediate if the
second assertion of~(\ref{dichotomy}) always holds, but there is a
priori no reason for the first assertion not to hold. Remember that $\min_{[-R,R]}v_n(t^*_n,\cdot)>0$ for each fixed $n\ge N$. So, if (\ref{claimrho}) were not true,
then, up to extraction of a subsequence, there would exist a sequence
of points $(y_n)_{n\ge N}$ in $[-R,R]$ such that
\[
v_n(t^*_n,y_n)\to0\hbox{ and }y_n\to y_{\infty}\in[-R,R]\hbox{ as
}n\to+\infty.
\]
We will use the translates 
\[
w_n(t,x)=v_n(t+t^*_n,x),
\]
defined for all $n\ge N$ and $(t,x)\in(-t_n-t^*_n,+\infty)\times\R$.
Since the functions $v_n$   solve
\begin{equation}\label{eqvn}
(v_n)_t=(v_n)_{xx}+c_n(v_n)_x+\mu\,v_n\,(1-\phi*v_n)\ \hbox{ in }(-t_n,+\infty)\times\R,
\end{equation}
the functions $w_n$ solve the same equation, in
$(-t_n-t^*_n,+\infty)\times\R$. Notice that $-t_n-t^*_n\le-1$ for
all~$n\ge N$, that $c_n\to c_{\infty}$ as $n\to+\infty$, that the
functions $w_n$ are all nonnegative and that the
sequence~$(\|w_n\|_{L^{\infty}((-t_n-t^*_n,+\infty)\times\R)})_{n\ge
  N}$ is bounded. Therefore, from standard parabolic estimates, the
functions $w_n$ converge in $C^{1,2}_{loc}((-1,+\infty)\times\R)$, up
to extraction of a subsequence, to a classical bounded solution
$w_{\infty}$ of
\[
(w_{\infty})_t=(w_{\infty})_{xx}+c_{\infty}(w_{\infty})_x
+\mu\,w_{\infty}\,(1-\phi*w_{\infty})\ \hbox{ in
}(-1,+\infty)\times\R,
\]
such that 
\[
\hbox{$w_{\infty}(t,x)\ge0$ for all $(t,x)\in(-1,+\infty)\times\R$
and $w_{\infty}(0,y_{\infty})=0$.}
\]
It follows form the strong maximum
principle and the uniqueness of the solutions of the Cauchy problem
that
\[
\hbox{$w_{\infty}(t,x)=0$ for all $(t,x)\in(-1,+\infty)\times\R$.}
\]
In
other words, $w_n\to0$ as $n\to+\infty$ (at least) locally uniformly
in~$(-1,+\infty)\times\R$, whence 
\[
\hbox{$\phi*w_n\to0$ as $n\to+\infty$,}
\]
locally uniformly in $(-1,+\infty)\times\R$ by boundedness of the
sequence $(\|w_n\|_{L^{\infty}((-1,+\infty)\times\R)})_{n\ge
  N}$. Therefore, we have
\[
\hbox{$v_n(t^*_n,\cdot)\to0$ and
$(\phi*v_n)(t^*_n,\cdot)\to0$ locally uniformly in $\R$ as
$n\to+\infty$.}
\]
This is a contradiction to (\ref{dichotomy}), since
both $\delta$ and $\eta$ are positive. Therefore,
claim~(\ref{claimrho}) is proved.

Now, (\ref{delta2}),~(\ref{claimrho}) and~(\ref{eqvn}) imply that we
have the following situation: for every $n\ge N$, one has $-t_n+1\le t^*_n<0$
and, in the box $[t^*_n,0]\times[-R,R]$, the nonnegative
function~$v_n$ satisfies 
\begin{equation}\label{eqvnbis}\left\{\baa{rcll}
  (v_n)_t & = & (v_n)_{xx}+c_n(v_n)_x+\mu\,v_n\,(1-\phi*v_n)\vspace{3pt}\\
  & \ge & (v_n)_{xx}+c_n(v_n)_x+\mu\,(1-\delta)\,v_n & \hbox{in }[t^*_n,0]\times[-R,R]\vspace{3pt}\\
  v_n(t,\pm R) & \ge & 0 & \hbox{for all }t\in[t^*_n,0],\vspace{3pt}\\
  v_n(t^*_n,x) & \ge & \rho & \hbox{for all }x\in[-R,R].\eaa\right.
\end{equation}
On the other hand, for every $n\ge N$, the  function
$\psi_n$ defined in $[-R,R]$ by
\[
\psi_n(x)=\rho\,e^{-c_nx/2-cR/2}\cos\Big(\frac{\pi x}{2R}\Big)
\]
satisfies $0\le\psi_n\le\rho$ in $[-R,R]$ from~(\ref{defcn}), $\psi_n(\pm R)=0$ and
\[
\psi_n''+c_n\psi_n'+\mu\,(1-\delta)\,\psi_n=
\Big(\mu\,(1-\delta)-\frac{c_n^2}{4}-\frac{\pi^2}{4R^2}\Big)\psi_n\ge0\
\hbox{ in }[-R,R]
\]
from our choice of $\delta$ and $R$: see
(\ref{defcn}),~(\ref{defdelta}) and~(\ref{defR}). In other words, the
time-independent function $\phi_n$ is a subsolution for the
problem~(\ref{eqvnbis}) satisfied by $v_n$ in
$[t_n^*,0]\times[-R,R]$. 
It follows then from the parabolic maximum
principle that
\[
v_n(t,x)\ge\phi_n(x)\ \hbox{ for all }(t,x)\in[t^*_n,0]\times[-R,R],
\]
and for all $n\ge N$. In particular,
\[
u(t_n,x_n)=v_n(0,0)\ge\phi_n(0)=\rho\,e^{-cR/2}\ \hbox{ for all }n\ge
N.
\]
However, assumption~(\ref{utnxn}) means that $u(t_n,x_n)\to0$ as
$n\to+\infty$, while $\rho\,e^{-cR/2}>0$ from the positivity of $\rho$
in~(\ref{claimrho}). One has then reached a contradiction.

Therefore,~(\ref{spreadinf}) holds and the proof of
Theorem~\ref{th3} is complete.\hfill$\Box$


\end{document}